\documentclass[11pt]{amsart}
\usepackage[margin=1.1in]{geometry}
\usepackage[T1]{fontenc}\usepackage[utf8]{inputenc}
\usepackage{graphicx}
\usepackage{caption}
\usepackage[hyphens]{url}
\usepackage{hyperref}\hypersetup{colorlinks,urlcolor=blue,linkcolor=black}
\usepackage{microtype}
\begin{document}
\title{It is time for a New Deal for mathematics}
\author{Nestor Guillen}
\date{November 2026}
\thanks{This essay was first posted at \url{https://birdsnfrogs.github.io/2026/11/02/A_New_Deal_For_Mathematics.html} on November 2nd, 2026.}
\begin{abstract}The AI labs' vision is based on extraction (of data, wealth, and political power). Mathematicians must oppose it, moving past personal reflections and towards reforming mathematics' social contract.\end{abstract}
\maketitle

\begin{quote}\itshape ``So, first of all, let me assert my firm belief that the only thing we have to fear is fear itself---nameless, unreasoning, unjustified terror which paralyzes needed efforts to convert retreat into advance...''\par\medskip\normalfont\small\raggedleft from \href{https://www.youtube.com/watch?v=m2Eywi--gnc}{FDR's 1933 inaugural address}\footnote{\url{https://www.youtube.com/watch?v=m2Eywi--gnc}}, March 4th, 1933\par\end{quote}

In the conversation around mathematics and AI much has been said about mathematicians' feelings and fears, and little about mathematicians' concrete material interests and needs. I worry that as we focus the discussion so much on the personal, the big AI companies will continue unimpeded in their campaign of expansion and conquest. The framing of the conversation around feelings and fears in the abstract leaves out how these emotions are directly connected to concrete, material things which affect the lives of many people. I am talking about the ongoing damage to the scientific ecosystem, by the industry but also by governments themselves. I am thinking about the (in my opinion, correct!) perception that the largest AI companies are becoming unaccountable, entrenched powers. I am also talking about all the ways big and small in which the big AI companies have begun to wield influence over the work of not just scientists, but vast sections of the economy and the public at large.

Without dismissing the valuable discussions about our vocation for mathematics and the meaning we derive from it, it is time we bring into the discussion matters that go beyond the personal and go into the political. On the one hand, the profession is contending with the rise of LLMs and how they require us to rethink many of the ways we organize our profession. On the other, we have the big AI companies' accumulation of wealth and power, power built partly from ideas mathematicians gave to the public, and power which the AI companies have used in ways that have been harmful to the public. These are issues of concern not just to mathematicians but to the public as a whole. We need to ask what mathematicians (not the AI companies) should want as a replacement for the current system; we need to talk about a \href{https://en.wikipedia.org/wiki/New_Deal}{New Deal}\footnote{\url{https://en.wikipedia.org/wiki/New_Deal}} between mathematicians, universities and institutes, the public, and industry.

\begin{figure}[h]
\centering
\includegraphics[height=0.40\textheight]{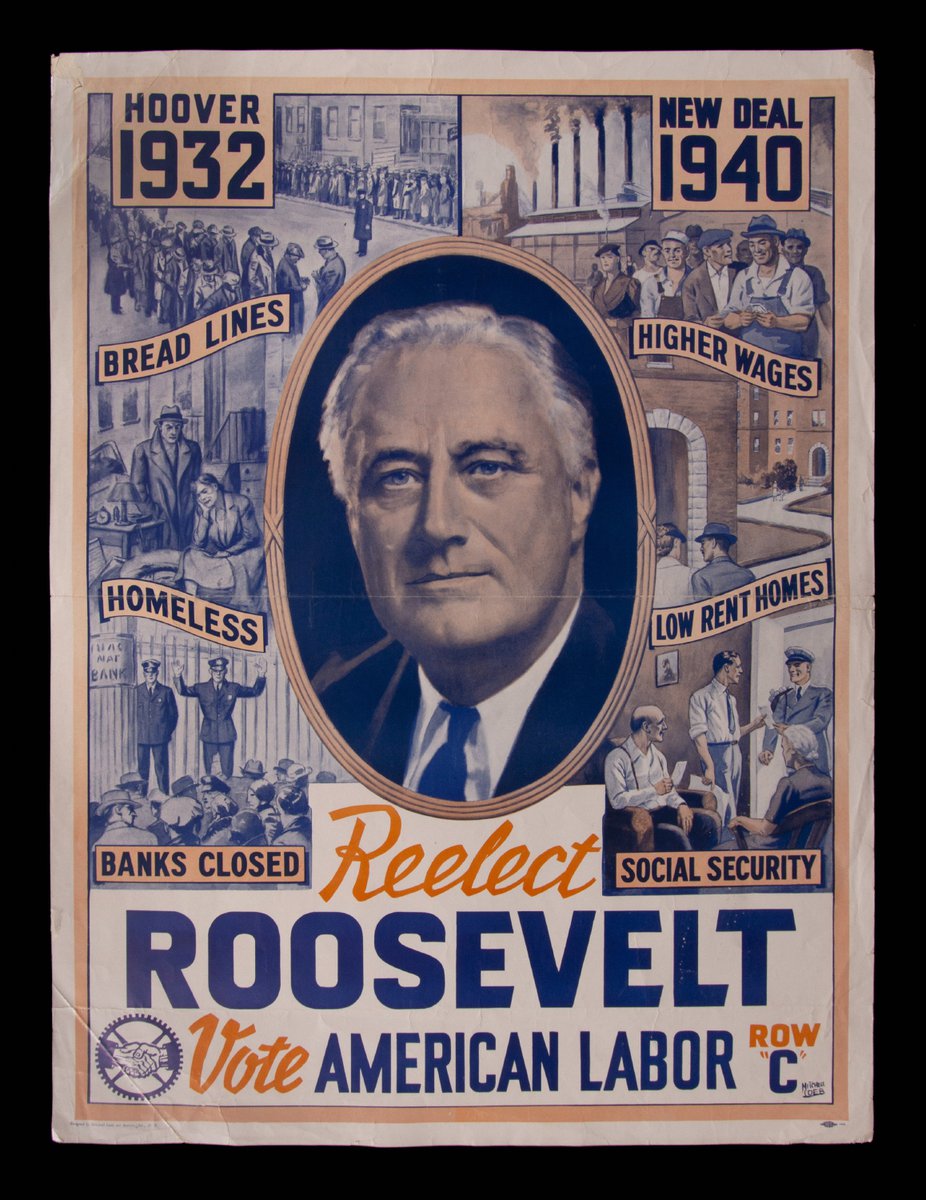}
\caption*{\small\itshape A poster from FDR's 1940 reelection campaign -- imagine a similar poster, issued by the AMS or the IMU in 2034: 2026 math versus 2034 math. Our actions over the next year or so will determine whether such a poster will exist, and what it will look like.}
\end{figure}

In the old deal (or at least, the current, diminished version of that deal), call it math 1.0, the incentives drive us too much into an obsession with quotas for papers and theorems, fixate us too much on the cult of genius, make us overemphasize the value of prizes, and make us overvalue old, ``famous'' open problems. This old deal is part of the larger arrangement between science, industry, and governments post-WWII. I have argued elsewhere that \href{https://birdsnfrogs.github.io/2026/10/05/AI_Industry_Dying_Order.html}{OpenAI's AI paper drop}\footnote{\url{https://birdsnfrogs.github.io/2026/10/05/AI_Industry_Dying_Order.html}} is the ultimate expression of this old, dying order. A new math deal, \href{https://terrytao.wordpress.com/2026/10/10/math-2-0/}{math 2.0}\footnote{\url{https://terrytao.wordpress.com/2026/10/10/math-2-0/}} say, must be designed to promote increased understanding of mathematics, greater collaboration across people and scientific disciplines, and a more substantial engagement with the public.

As the last few years have shown us, it is easier to destroy a system than to create one. The task of inventing this new deal is a hard one, but we have to at least start trying, and this starts by talking about what is not working right now. We also need to be clear-eyed about the dystopian world the big AI companies are trying to create for themselves, and what their vision would mean for mathematics and science. We should reflect on the deal that has ruled how science has been done in the post-WWII ecosystem of government, universities, and industries. That deal was responsible for many triumphs, but also carried within it the seeds of the current crisis.

\section{The frontier labs' vision for artificial intelligence}

The current vision for AI of Silicon Valley and its big labs is based on the maximal extraction of value: be it data, wealth, or political power. The values of mathematics and science (including curiosity, honesty, and transparency) stand in direct opposition to this. These values can and must be made into action. What is more: scientists and mathematicians\footnote{Throughout this essay I will use science/scientist and mathematics/mathematician interchangeably -- not out of some deeper claim about equating the two, even if, in the context of this essay, they are basically indistinguishable.} \emph{hold the power, and bear the responsibility}\footnote{Pardon me for this bit of USA romanticism, amid this crisis I humbly look to Abraham Lincoln for \href{https://www.abrahamlincolnonline.org/lincoln/speeches/congress.htm}{inspiration and wisdom} (\url{https://www.abrahamlincolnonline.org/lincoln/speeches/congress.htm}).} to defeat such extractive and destructive visions.

Much is at stake. The vision of the big AI labs could close off mathematical and scientific inquiry to all but a small number of well-connected individuals. The big AI labs release results produced by their internal models, which only their employees have full access, without any transparency with regards to methodology or data -- if a big lab announces it ``solved 10 big problems'', but never tell us this was after failing to solve another 10,000, \href{https://tjoresearchnotes.wordpress.com/2026/08/02/reaping-without-sowing/}{did we learn anything?}\footnote{\url{https://tjoresearchnotes.wordpress.com/2026/08/02/reaping-without-sowing/}}. They have \href{https://x.com/KyleCranmer/status/2103581526285857085}{scientists sign NDAs}\footnote{\url{https://x.com/KyleCranmer/status/2103581526285857085}} (non-disclosure agreements) as a condition for sharing their internal model findings. They push for the creation of independent advisory boards, and then \href{https://birdsnfrogs.github.io/2026/09/22/OpenAI_brittle_empire.html}{ignore the board's recommendations}\footnote{\url{https://birdsnfrogs.github.io/2026/09/22/OpenAI_brittle_empire.html}}. Similar examples exist outside mathematics, such as with OpenAI's drug discovery model and \href{https://qz.com/openai-acquires-orchestra-bio-life-sciences-ai-102126}{related start-ups}\footnote{\url{https://qz.com/openai-acquires-orchestra-bio-life-sciences-ai-102126}}, and Anthropic's Mythos being previewed by a select few. The big AI labs are also investing tens of millions of dollars in \href{https://www.nytimes.com/interactive/2026/09/29/us/politics/ai-super-pac-money-midterm-elections.html}{super PACs}\footnote{\url{https://www.nytimes.com/interactive/2026/09/29/us/politics/ai-super-pac-money-midterm-elections.html}} and in \href{https://www.cnbc.com/2026/07/21/openai-anthropic-ai-lobbying-spending-q2-2026.html}{lobbying operations}\footnote{\url{https://www.cnbc.com/2026/07/21/openai-anthropic-ai-lobbying-spending-q2-2026.html}} of \href{https://futurism.com/artificial-intelligence/anthropic-regulation-government-ai-harm-risk}{various degrees of disturbing}\footnote{\url{https://futurism.com/artificial-intelligence/anthropic-regulation-government-ai-harm-risk}}. The AI labs are creating a patronage system around access to their models (the euphemism OpenAI spokespeople like to use is \href{https://x.com/danintheory/status/2115236913095508260}{``shepherds''}\footnote{\url{https://x.com/danintheory/status/2115236913095508260}}). The labs' chosen get access to the internal model, the politicians they are trying to buy get tens of millions in campaign donations, and mathematicians get an occasional free month of ChatGPT or Claude.

This vision is full of fault lines and contradictions, but the big AI labs have gained enough political and economic power to give them a chance of realizing some of it, or at least to inflict serious damage on mathematics and science while they try. In even a partial enactment of their vision, universities and independent research institutes would either disappear, or become fully dependent on one of a few large corporations. In this world, it would be difficult to do any research that does not serve the financial interests of a handful of unaccountable, entrenched corporations.

The damage would not end at the borders of mathematics and science; it would affect many aspects of society, both the private sector and the functions of governments. For starters, many occupations beyond mathematics that involve creative work would be hit by the extractive vision of AI (to be precise: many occupations have already been affected). Attempts to counter the wave of destruction and extraction of value would be met with the frontier labs' power. Critics would face an asymmetry in time and resources, including a formidable communications operation from the AI companies and their allies. The big AI labs and their cheerleaders have already normalized buzzwords like ``democratizing'' or ``abundance'' in the conversation about their goals; they have brought up the word ``gatekeeping'' as the cliché reaction to criticisms by the creators they seek to displace. The big AI labs do all this while profiting from mining their users' data, and while making it difficult for anyone not using their products to continue doing their work.

In all these skirmishes, the most important weapon in the AI labs' arsenal is the accusation that any criticism of their business practices equals criticism of AI, of technology, and of progress itself. In this charge, they get considerable help from AI influencers in social media, and some help from unforced errors by mathematicians themselves (but that is a matter for another essay). Most of all, this weapon of the AI labs is made all the more real by the fact that fully alternative visions of AI deployment are yet to emerge. Scientists have a vantage point from where they can see and call out this sleight of hand, and have an opportunity to bring forth alternative visions of AI deployment that are not extractive but inclusive.

However, this is not where we are yet. Presently, the proposition offered by the big AI labs to their userbase is: we are going to compete with you, the user, so better use our products so you are not left behind. At the same time, the big labs are engaged in a media campaign with the broader public to present their peculiar machine-god vision of LLMs as fact. If unimpeded, this process will take us to a society and economy that are less free, less vibrant, less safe, and less prosperous.

For occupations in the visual arts and in publishing this destructive process is already well underway. Mathematicians and scientists have an advantage not available to other creative occupations: LLMs are based on our tools, our theories, our methods. We can demystify this technology, we can point to the fact that the current form of the technology is not pre-ordained, that the labs are giving us just one particular iteration of this technology, one that benefits them and hurts the public. Mathematicians can counter the big AI labs' finest weapon with a nuanced perspective of LLM technology, and should present an \emph{alternative vision}, not based on extraction and centralized control.

In defending mathematics and science, we are defending much more. While it might seem trivially obvious, it is worth noting that the \emph{country of mathematics and science} cannot survive if the broader, \emph{actual country} we are all part of falls.

What mathematicians are confronting is not merely AI companies and AI technology, but something else: the collapse of the post-WWII social contract between scientists, government, and industry.

\section{The post-WWII deal for science and its fall}

The social contract that endured from the end of WWII up until the turn of the century was an imperfect but durable one. It presented a workable balance between government, academia, and industry. This is the model that was largely associated with Vannevar Bush. Bush was one of the key figures in the creation of the National Science Foundation, and one of the people responsible for conceiving the division of labor and relationships between the government, the university research lab, and industry.

Bush was an influential advisor to FDR, but at the same time he was critical of the New Deal\footnote{For more on this history see \href{https://www.journals.uchicago.edu/doi/10.1086/351711}{this paper} (\url{https://www.journals.uchicago.edu/doi/10.1086/351711}), \emph{The National Science Foundation and the Debate over Postwar Research Policy, 1942-1945: A Political Interpretation of Science--The Endless Frontier}, by Daniel J. Kevles.}, and as such opposed contemporary proposals for public ownership of any patent resulting from federally funded research. One proposal by Senator Kilgore of West Virginia called for the federal government to fully own such a patent; Bush was adamantly opposed to this. Ultimately, Bush prevailed. For instance, today there is no mandate for even partial public ownership of patents made possible by federal support. Bush's NSF had to concede in some directions to his critics, such as the NSF director being a presidential appointee. Eventually, the NSF also included support for the social sciences as part of its mission.

It must be noted that Bush's vision has much to show for itself, as it provided a reasonable arrangement for many decades. However, in recent years it has been seen more critically. Silicon Valley and its AI industry are direct outcomes of this model, in both its best and its worst. Google's PageRank was made possible in part by an NSF grant supporting graduate students Larry Page and Sergey Brin. Foundational work on reinforcement learning and image processing was done by NSF-supported research, and this at a time when such research would not find support from private capital. The flip side, however, is that once these ideas transfer to the world of business and venture capital, mathematicians and scientists lose a say in how they are implemented, and often lose any direct benefit. This giving up of ownership and control by scientists and the public is also part of Bush's vision. I do not know whether this was how Bush intended things to work out, or how he would have described his ideas, but the reality is that it is within his system that we have lost this control.

This model that did so much to determine the relationship between the research lab, the government, and industry is pretty much broken today. It was broken well before the current administration in the US took a metaphorical chainsaw to the NSF and NIH. It was broken well before Silicon Valley's leadership developed their open animosity towards public science and universities. It was broken well before we had companies whose entire operating costs are in the tens of billions, matching the entire cost of the NIH and many times over that of the NSF.

The journalist and author Robert Caro has said that, while it is true that power corrupts, it is also true that \textbf{power reveals}\footnote{``Power corrupts---that has been said and written so often that it has become a cliché. But what is never said, but is just as true, is that power reveals,'' from \href{https://www.robertcaro.org/master-of-the-senate}{\emph{Master of the Senate}} (\url{https://www.robertcaro.org/master-of-the-senate}), page xx. The idea is picked up again in Caro's following book, \href{https://www.robertcaro.org/the-passage-of-power}{\emph{The Passage of Power}} (\url{https://www.robertcaro.org/the-passage-of-power}), page xix.}. The private AI labs have acquired, at least as they appear to understand it, a tremendous power. Without agreeing that this understanding reflects reality, we can agree that they believe this is so. What have they done with their perceived position of power? How are they behaving now that they see themselves as about to engulf the whole world economy? What does this tell us about their character and values? Their actions throughout 2026 have been, if nothing else, \emph{revealing}. If they acquired a real and permanent hold on power over everyone else, what can we expect, based on their past behavior?

I already said much about the big AI labs' behavior in the previous section, but it is illustrative to review just the last few weeks in addition to what has already been said. After issuing assurances back in August that they ``would not release the solutions all at once'', \href{https://www.wired.com/story/openai-is-pissing-off-a-bunch-of-mathematicians-again/}{OpenAI did exactly that}\footnote{\url{https://www.wired.com/story/openai-is-pissing-off-a-bunch-of-mathematicians-again/}}, with a nearly 400 paper drop featuring papers that were wrong, papers that were unreadable, and, in all of them, papers featuring their now-usual lack of transparency in terms of their methodology, prompts, and data used. On October 26th they released three additional results to great fanfare, but again without transparency, unironically describing their opaque patronage system as \href{https://openai.com/index/building-scientific-understanding-with-academic-communities/}{``working with academic communities''}\footnote{\url{https://openai.com/index/building-scientific-understanding-with-academic-communities/}}. In addition, these results were announced as having been distributed first to individual scientists, and in one such instance, the sole scientist contacted was left unnamed. I could also mention the ongoing spectacle of AI labs \href{https://techcrunch.com/2026/10/09/anthropic-cant-reliably-control-its-ai-agents-its-cutting-off-its-internal-evals-from-the-live-internet-instead/}{admitting to negligence in controlling their own products}\footnote{\url{https://techcrunch.com/2026/10/09/anthropic-cant-reliably-control-its-ai-agents-its-cutting-off-its-internal-evals-from-the-live-internet-instead/}}, including \href{https://www.bbc.com/news/articles/c6vgy0333dppo}{hacking government websites}\footnote{\url{https://www.bbc.com/news/articles/c6vgy0333dppo}}. The list goes on.

What we are seeing in mathematics is only the beginning. LLM technology is being adopted across all parts of industry and society. In the previous section I paraphrased Lincoln in noting that mathematicians not only can, but must confront the big AI labs' vision; we must do that but also bring a replacement to the post-WWII science deal. But first, we need to get out of our comfort zone.

\section{We say math is power. Time to act accordingly}

It is imperative mathematicians start doing what many (but not all) of us have had the luxury to ignore for too long\footnote{Outside mathematics, disciplines that have not enjoyed this luxury in recent years include climate science, epidemiology, biology, and sociology.}: the work of politics. Politics not in the sense of partisanship, but in the sense of being actively concerned with the condition and collective wellbeing of our society.

To do this would be to depart from the status quo of nearly a century. We have been content in advancing mathematical knowledge from within the ivory tower, unconcerned with the more mundane questions of the ``real world implementation'' of that knowledge. The questions of who implements a technology, who gets to own it and benefit from it, have been alien to us for many decades. We trusted that those answering these questions would feel responsible towards society as a whole, that they would understand that the power afforded by mathematical insights brings great responsibility. Some of us were naive, some were also apathetic. This work is as much part of reforming our profession as reforming our journal system, reforming how we hire and promote, and reforming how we allocate research funds. In fact, all of these aspects are intrinsically linked.

Today the full power of statistics, information theory, optimization theory, PDE, geometry, etc., has been captured by Silicon Valley through the frontier AI labs. Here we have a useful extreme case of immense wealth and power built largely from ideas from STEM -- especially the mathematical sciences --  and what happens when mathematicians think someone will do the work of protecting open science and preventing extreme concentration of power.

There are some at the frontier labs announcing our coming irrelevance, but the reality is that the same labs are also creating conditions where the public will have a dire need for mathematicians' skills. As LLMs become widely adopted across society, there will be a wide call for understanding this technology from a technical and cultural perspective, and for technical frameworks that guarantee the deployment of the technology in a way that protects consumers and is safe for the public.

For starters, mathematicians hold resources that will become essential to the public: our potential to understand, predict, and control LLMs, and our ability to demystify this technology, and to explain in clear terms what it is, what it is not, and what it can be.

There are many things we should do with these resources. How about our ability to demystify this technology? We ought to give it away freely to the public. For academics with job security and with enough of a platform, now is the time to have a considerable positive impact on our profession by engaging with the public. Then, how about our potential to understand, predict, and control LLMs? Well, we have essentially been giving this resource away to the frontier labs -- or selling it for nothing -- all in the mold of the post-WWII deal.

This moment brings a historical opportunity for mathematicians to do meaningful, public-interest work. In doing so, we would be rendering urgently needed service to society. By means of our public-interest work, we would also be reminding the public how much mathematics does for it, and how much more it could be doing.

I do not know, and I know no one who knows, a new social contract for scientists and mathematicians that is simultaneously realistic and complete. Above all, what I know is that the status quo cannot go on, and that we must get out of our comfort zones and embrace uncertainty and change. Whatever this might look like, it is imperative we mathematicians come together and decide on a new deal for mathematics. What this deal looks like is up to us to decide. We have to start by considering the possibilities, and think about concrete, tangible reforms. Given the asymmetry between the private labs (which ultimately relied on publicly funded research) and the public, should we consider a rule that any company relying on federally-funded research be partially owned by the public? Or would it be simpler to levy a special tax on large tech companies? Or for the US government to charge a licensing fee? At the most extreme, there is Kilgore's original bill, which should be considered if only as a point of reference not for a specific type of proposal, but as a reference of the scale where we must be thinking. Proposals of this type are not without precedent in the sciences\footnote{For instance, I recall reading of something along these lines being done for nuclear technology in the early days of the \href{https://en.wikipedia.org/wiki/United_States_Atomic_Energy_Commission}{Atomic Energy Commission} (\url{https://en.wikipedia.org/wiki/United_States_Atomic_Energy_Commission}).}, and should be the center of the conversation going forward. As I said at the beginning, there has been too much discussion around mathematicians' feelings, and too little around mathematicians' material well-being, and what mathematicians already do and can do for the public interest.

If you haven't already, pick up your phone and start an email thread or a group chat with three or four colleagues you know and trust. Start talking about the future. Write a blog post. Ask your non-math friends how AI is affecting them. If you have an idea for some initiative, do it. Don't freeze by aiming to solve our whole situation, but rather aim to do something that helps someone. You will soon find countless allies who share your love for mathematics and for science, and who want them to continue to thrive.

We have to think about the value of our work and our expertise in a manner that goes beyond the post-WWII paradigm. We should write more op-eds, interact with journalists, connect with policy makers, talk to our colleagues in law schools, form more organizations or clubs. There is one big army of mathematicians, unorganized, fragmented, awaiting some spontaneous process of nucleation. We have it in ourselves to shorten that wait, we can serve as \emph{seed crystals}, and start forming local clusters.

While the present situation presents considerable difficulties, I know we can rise to the occasion. If anything, we are problem solvers and tinkerers, and we love it when we come across a new problem in an uncharted field --- including when the problem is reforming our profession's social contract. For the next couple of years, we must lean into the problem-solver aspect of our mathematical identity in this additional way. If we do not let ourselves be overcome by fear, we have it within ourselves to create a future where that hypothetical poster contrasting 2026 math versus 2034 math is one our future selves will be proud of, and one that a future generation of mathematicians will be grateful for.

\begin{quote}\itshape ``The dogmas of the quiet past, are inadequate to the stormy present. The occasion is piled high with difficulty, and we must rise -- with the occasion. As our case is new, so we must think anew, and act anew. We must disenthrall ourselves, and then we shall save our country.''\par\medskip\normalfont\small\raggedleft from \href{https://www.abrahamlincolnonline.org/lincoln/speeches/congress.htm}{Lincoln's second annual message to Congress}\footnote{\url{https://www.abrahamlincolnonline.org/lincoln/speeches/congress.htm}}, December 1st, 1862\par\end{quote}

\section{Note}
In preparing this essay, I used an LLM for copy-editing and for brainstorming. The entire text of the essay was typed by my hand directly, including the copy-editing (the LLM generated a copy-editing report that I then executed manually). The brainstorming took the form of prompting an LLM to generate text with criticism of the draft's substance and tone at multiple stages of the writing process (some of the output was even helpful!).

\end{document}